# Rapidly Convergent Series of the Divisors Functions

N. A. Carella

*Abstract*: This note gives a few rapidly convergent series representations of the sums of divisors functions. These series have various applications such as exact evaluations of some power series, computing estimates and proving the existence results of some special values of the sums of divisors functions.

## 1. Introduction

Let $N \in \mathbb{N}$ be an integer, and let $s \in \mathbb{C}$ be a complex number. The complex valued divisors function is defined by $\sigma_s(N) = \sum_{d \mid N} d^s$. It has a power series expansion (Ramanujan series) given by

$$\frac{\sigma_{s-1}(N)}{N^{s-1}} = \zeta(s) \sum_{k=1}^{\infty} \frac{c_k(N)}{k^s}, \tag{1}$$

where $c_k(N) = \sum_{\gcd(u,k)=1} e^{-i2\pi u N/k}$, see [HD], [SR]. The representatation of the the zeta function used in (1) depends on the domain of definition:

$$\zeta(s) = \begin{cases} \sum_{n \geq 1} n^{-s} & \text{if } \Re e(s) > 1, \\ \left(1 - 2^{1-s}\right)^{-1} \sum_{n \geq 1} (-1)^{n-1} n^{-s} & \text{if } \Re e(s) > 0, \end{cases} \tag{2}$$

see [TT]. The special cases for $s \in \mathbb{N}$ are the finite sums. For example, the first three sums of divisors functions are

$$\frac{\sigma(N)}{N} = \sum_{d \mid N} \frac{1}{d}, \qquad \frac{\sigma_2(N)}{N^2} = \sum_{d \mid N} \frac{1}{d^2}, \qquad \frac{\sigma_3(N)}{N^3} = \sum_{d \mid N} \frac{1}{d^3}, \tag{3}$$

respectively. The goal of this note is to develop a few rapidly convergent series expansions of the sums of divisors functions. Some of the series expansions are the followings.

**Theorem 1.** Let $N \geq 1$ be an integer, and let $\sigma_s(N)$ be a sum of divisors function. Then

(i) $\dfrac{\sigma(N)}{N} = 3 \sum_{n=1}^{\infty} \dfrac{1}{n^2} \sum_{d \mid n} \binom{2d}{d}^{-1} c_{n/d}(N),$

(ii) $\dfrac{\sigma_2(N)}{N^2} = \dfrac{5}{2} \sum_{n=1}^{\infty} \dfrac{1}{n^3} \left[ \sum_{d \mid n} (-1)^{d-1} \binom{2d}{d}^{-1} c_{n/d}(N) \right],$ \hfill (4)

(iii) $\dfrac{\sigma_3(N)}{N^3} = \dfrac{36}{17} \sum_{n=1}^{\infty} \dfrac{1}{n^4} \sum_{d \mid n} \binom{2d}{d}^{-1} c_{n/d}(N),$

where the symbol $\binom{n}{k} = \dfrac{n!}{k!(n-k)!}$ denotes the $n$th binomial coefficient.

**Proof:** These series spring from the convolution of the rapidly convergent series of the zeta function in Theorem 2, and the Ramanujan series in (1). That is,

$$\dfrac{\sigma_{s-1}(N)}{N^{s-1}} = \zeta(s) \cdot \sum_{n=1}^{\infty} \dfrac{c_n(N)}{n^s}. \tag{5}$$

A routine calculation produces the claim. ∎

The $n$th coefficient

$$a_n = \sum_{d \mid n} \binom{2d}{d}^{-1} c_{n/d}(N) \quad \text{or} \quad a_n = \sum_{d \mid n} (-1)^{d-1} \binom{2d}{d}^{-1} c_{n/d}(N), \tag{6}$$

is a weighted sieving gadget. It sieves the factors of $N \geq 1$ up to a weight depending on the parameter $\gcd(n, N)$. Its maximal value, and its minimal value are

$$a_n = \sum_{d \mid n} \binom{2d}{d}^{-1} \varphi(n/d) \quad \text{and} \quad a_n = \sum_{d \mid n} \binom{2d}{d}^{-1} \mu(n/d), \tag{7}$$

if $\gcd(n, N) = n$ or if $\gcd(n, N) = 1$ respectively. This follows from the identity $c_k(N) = \mu(k/e) \dfrac{\varphi(k)}{\varphi(k/e)}$, $e = \gcd(k, N)$, see [AP, p. 162]. Here $\mu(n) \in \{-1, 0, 1\}$ is the Mobius function, and $\varphi(n) = n \prod_{p \mid n} (1 - 1/p)$ is the totient function. This identity expresses $c_k(N) \in \mathbb{Z}$ as an integer.

The first power series expansion of the sum of divisors function in Theorem 1 is expressed in terms of rational numbers:

$$\frac{\sigma(N)}{3N} = \frac{1}{2} + \frac{\binom{2}{1}^{-1} c_2 + \binom{4}{2}^{-1}}{4} + \frac{\binom{2}{1}^{-1} c_3 + \binom{6}{3}^{-1}}{9} + \frac{1}{16}\left(\binom{2}{1}^{-1} c_4 + \binom{4}{2}^{-1} c_2 + \binom{8}{4}^{-1}\right) + \cdots, \tag{8}$$

using $c_1 = 1$, and $c_k = c_k(N)$. It has significantly faster rate of convergence than the earlier Ramanujan series

$$\frac{\sigma(N)}{N} = \frac{\pi^2}{6}\left(1 + \frac{c_2(N)}{4} + \frac{c_3(N)}{9} + \frac{c_4(N)}{16} + \cdots\right). \tag{9}$$

## 2. Rapidly Convergent Series of the Zeta Function

The zeta function is the complex valued function defined by

$$\zeta(s) = \sum_{n \geq 1} \frac{1}{n^s} \tag{10}$$

on the complex half plane $\{ s \in \mathbb{C} \colon \mathfrak{Re}(s) > 1 \}$. It has a wide range of different rapidly convergent series.

***Theorem 2.*** Let $s = 2, 3$ or $4$. Then

(i) $\zeta(2) = 3 \sum_{n=1}^{\infty} \binom{2n}{n}^{-1} \frac{1}{n^2}$,

(ii) $\zeta(3) = \frac{5}{2} \sum_{n=1}^{\infty} (-1)^{n-1} \binom{2n}{n}^{-1} \frac{1}{n^3}$, \hfill (11)

(iii) $\zeta(4) = \frac{36}{17} \sum_{n=1}^{\infty} \binom{2n}{n}^{-1} \frac{1}{n^4}$.

***Proof***: See [AT], [BB], [BR], [LR], [DL, 25.6.7], et alii. ∎

For $s > 4$, this family of rapidly convergent series of the zeta function is expected to be of the form

$$\text{(i)} \ \zeta(s) = c_s \sum_{n=1}^{\infty} \binom{2n}{n}^{-1} \frac{1}{n^s}, \quad \text{or} \quad \text{(ii)} \ \zeta(s) = c_s \sum_{n=1}^{\infty} (-1)^{n-1} \binom{2n}{n}^{-1} \frac{1}{n^s}, \tag{12}$$

where $c_s \in \mathbb{Q}$ is a rational constant. But the constants $c_s > 0$ are unknown, see [AM], [BB], [BR]. Other families of slightly more complex series are generated by the generating functions

(i) $\sum_{k=0}^{\infty} \zeta(2k+3) z^{2k} = \sum_{n=1}^{\infty} \frac{1}{n^3(1-z^2/n^2)}$

$= \sum_{n=1}^{\infty} \binom{2n}{n}^{-1} \frac{(-1)^{n-1}}{n^3} \left( \frac{1}{2} - \frac{2}{1-z^2/n^2} \right) \prod_{1 \leq m < n} \left( 1 - \frac{z^2}{m^2} \right),$

and

(ii) $\sum_{k=0}^{\infty} \left( 1 - \frac{1}{2^k} \right) \zeta(2k+2) z^{2k} = \sum_{n=1}^{\infty} \frac{1}{n^2(1-z^2/n^2)}$ (14)

$= \sum_{n=1}^{\infty} \binom{2n}{n}^{-1} \frac{(-1)^{n-1}}{n^3} \left( -\frac{1}{2} + \frac{2}{1-z^2/n^2} \right) \prod_{1 \leq m < n} \left( 1 - \frac{z^2}{m^2} \right).$

These generating function were developed in [LR]. Several families of nonrational rapidly convergent series representations of the zeta function are known. For example, the interesting series

$$\zeta(s) = \frac{1}{2^s - 2} \sum_{n=1}^{\infty} \frac{(s+1)_{2n}}{(2n)!} \frac{\zeta(s+2n)}{2^{2n}},$$ (15)

where $(x)_n = x(x+1)\cdots(x+n-1)$ is the rising factorial, and $\Re(s) > 0$, [SH, p. 440], is a rapidly convergent irrational representation of the zeta function. There are many other families of rapidly convergent series of the zeta function, confer the literature for details.

## 3. Sparse and Rapidly Convergent Series of the Sums of Divisors Functions

Other series can also be utilized to construct rapidly convergent series of the sum of divisors functions. A few cases will be demonstrated below.

***Lemma* 3.** Let $N \geq 1$ be an integer, and let $\sigma(N) = \sum_{d|N} d$ be the sum of divisors function. Then

$$\left( 1 - \frac{6 \log(2)^2}{\pi^2} \right) \frac{\sigma(N)}{2N} = \sum_{n=1}^{\infty} \frac{1}{n^2} \sum_{d|n} 2^{-d} c_{n/d}(N).$$ (16)

***Proof*:** The identity $L_s(z) + L_s(1-z) = \zeta(2) - \ln(z)\ln(1-z)$, $0 < z < 1$, of the polylogarithm function $L_s(z) = \sum_{n=1}^{\infty} z^n n^{-s}$ evaluated at $s = 2$, $z = 1/2$, has the well known value

$$\sum_{n=1}^{\infty} \frac{1}{2^n n^2} = \frac{\zeta(2)}{2} - \frac{\log(2)^2}{2},$$ (17)

see [DL, 25.12.6]. Taking the Dirichlet convolution of the two power series (1) and (17) yields the claim. ∎

The next series expansion has the potential of being a sparse series, and possibly a lacunary expansion of the sum of divisors function depending on the integer arguments $N \in \mathbb{N}$.

**Lemma 4.** Let $N \geq 1$ be an integer, and let $\sigma(N) = \sum_{d \mid N} d$ be the sum of divisors function. Then

$$\frac{90}{\pi^4} \frac{\sigma(N)}{N} = \sum_{n=1}^{\infty} \frac{1}{n^2} \sum_{d \mid n} \mu(d)^2 \, c_{n/d}(N). \tag{18}$$

**Proof**: Consider the generating series of the squarefree integers $\zeta(s)/\zeta(2s) = \sum_{n=1}^{\infty} \mu(n)^2 \, n^{-s}$. Taking the Dirichlet convolution of the power series (1) and this series yields

$$\frac{\zeta(s)}{\zeta(2s)} \frac{\sigma_{s-1}(N)}{\zeta(s) N^{s-1}} = \sum_{n=1}^{\infty} \frac{\mu(d)^2}{n^s} \sum_{n=1}^{\infty} \frac{c_n(N)}{n^s}. \tag{19}$$

A routine calculation produces the claim. ∎

## 4. Nonlinear Identities of the Sums of Divisors Functions

The different powers of the sums of divisors functions $\sigma_s(N) = \sum_{d \mid N} d^s$ can be identified with the simpler sum of divisors function $\sigma(N) = \sum_{d \mid N} d$ up to a weighted product of the prime divisors of the argument $N$.

**Lemma 5.** Let $N \geq 1$ be an integer, and let $\sigma_s(N) = \sum_{d \mid N} d^s$ be the sum of divisors function. Then

$$\frac{\sigma_s(N)}{N^s} = \frac{\sigma(N)}{N} \prod_{p^\alpha \| N} \left( \frac{p^{(s-1)(\alpha+1)} + p^{(s-2)(\alpha+1)} + \cdots + p^{\alpha+1} + 1}{p^{(s-1)\alpha+s-1} + p^{(s-1)\alpha+s-2} + \cdots + p^{(s-1)\alpha+1} + p^{(s-1)\alpha}} \right) \tag{20}$$

where $p^\alpha \| N$ denotes the largest prime power divisor of $N$.

**Proof**: Use the multiplicative formula to obtain

$$\frac{\sigma_s(N)}{N^s} = \frac{1}{N^s} \prod_{p^\alpha \| N} \left( \frac{p^{s(\alpha+1)} - 1}{p^s - 1} \right)$$

$$= \frac{1}{N^s} \prod_{p^\alpha \| N} \left( \frac{p^{\alpha+1} - 1}{p - 1} \right) \prod_{p^\alpha \| N} \left( \frac{p^{(s-1)(\alpha+1)} + p^{(s-2)(\alpha+1)} + \cdots + p^{\alpha+1} + 1}{p^{s-1} + p^{s-2} + \cdots + p + 1} \right)$$

$$= \frac{1}{N} \prod_{p^\alpha \| N} \left( \frac{p^{\alpha+1} - 1}{p - 1} \right) \frac{1}{N^{s-1}} \prod_{p^\alpha \| N} \left( \frac{p^{(s-1)(\alpha+1)} + p^{(s-2)(\alpha+1)} + \cdots + p^{\alpha+1} + 1}{p^{s-1} + p^{s-2} + \cdots + p + 1} \right)$$

$$= \frac{\sigma(N)}{N} \prod_{p^\alpha \| N} \left( \frac{p^{(s-1)(\alpha+1)} + p^{(s-2)(\alpha+1)} + \cdots + p^{\alpha+1} + 1}{p^{(s-1)\alpha + s - 1} + p^{(s-1)\alpha + s - 2} + \cdots + p^{(s-1)\alpha + 1} + p^{(s-1)\alpha}} \right).$$

(21)

These complete the proof. ∎

The two identities are

$$\frac{\sigma_2(N)}{N^2} = \frac{\sigma(N)}{N} \prod_{p^\alpha \| N} \left( \frac{p^{\alpha+1} + 1}{p^{\alpha+1} + p^\alpha} \right),$$

$$\frac{\sigma_3(N)}{N^3} = \frac{\sigma(N)}{N} \prod_{p^\alpha \| N} \left( \frac{p^{2(\alpha+1)} + p^{\alpha+1} + 1}{p^{2\alpha+2} + p^{2\alpha+1} + p^{2\alpha}} \right),$$

(22)

respectively. These identities can be used to develop the rapidly convergent series of the other divisors functions.

**Lemma 6.** Let $N \geq 1$ be an integer, and let $\sigma_s(N) = \sum_{d \mid N} d^s$ be the sum of divisors function. Then

(i) $\quad \dfrac{\sigma_2(N)}{N^2} = 3 \prod_{p^\alpha \| N} \left( \dfrac{p^{\alpha+1} + 1}{p^{\alpha+1} + p^\alpha} \right) \sum_{n=1}^\infty \dfrac{1}{n^2} \sum_{d \mid n} \binom{2d}{d}^{-1} c_{n/d}(N).$

(ii) $\quad \dfrac{\sigma_3(N)}{N^3} = 3 \prod_{p^\alpha \| N} \left( \dfrac{p^{2(\alpha+1)} + p^{\alpha+1} + 1}{p^{2\alpha+2} + p^{2\alpha+1} + p^{2\alpha}} \right) \times \sum_{n=1}^\infty \dfrac{1}{n^2} \sum_{d \mid n} \binom{2d}{d}^{-1} c_{n/d}(N).$

(23)